\documentclass{amsart}

\usepackage{amsmath}
\usepackage{amsfonts}
\usepackage{amssymb}
\usepackage{amsthm}
\usepackage{cite}

\usepackage{epsfig,epstopdf}
\usepackage[dvipsnames,usenames]{color}
\usepackage{hyperref}
\usepackage{xifthen}
\usepackage{graphicx}
\usepackage{tikz}
\tikzstyle{vertex}=[circle, draw, inner sep=0pt, minimum size=6pt]

\usepackage[
top    = 1in,
bottom = 1in,
left   = 1in,
right  = 1in]{geometry}
\usepackage{graphicx}

\tikzset{
  treenode/.style = {align=center, inner sep=0pt, text centered,
    font=\sffamily},
  arn_n/.style = {treenode, circle, black, text width=1em},
  arn_r/.style = {treenode, circle, blue, draw=blue, 
    text width=1em, very thick}}

\theoremstyle{plain}
\newtheorem{theorem}{Theorem}[subsection]
\newtheorem{lemma}[theorem]{Lemma}
\newtheorem{corollary}[theorem]{Corollary}
\newtheorem{proposition}[theorem]{Proposition}

\theoremstyle{definition}

\newtheorem{defi}[theorem]{Definition}
\newtheorem{example}[theorem]{Example}

\newcommand{\C}{\mathbb{C}}
\newcommand{\D}{\mathbb{D}}
\newcommand{\N}{\mathbb{N}}
\newcommand{\R}{\mathbb{R}}

\newcommand{\ds}{\displaystyle}

\newcommand{\calB}{\mathcal{B}}
\newcommand{\calF}{\mathcal{F}}

\newcommand{\dist}{\operatorname{d}}

\newcommand{\parent}[1]{\operatorname{par}({#1})}
\newcommand{\nparent}[2]{\operatorname{par}^{{#2}}({#1})}
\newcommand{\Chi}[1]{{\operatorname{Chi}}({#1})}
\newcommand{\nChi}[2]{{\operatorname{Chi}}^{#2}({#1})}

\newcommand{\conj}[1]{\overline{#1}}
\newcommand{\ceiling}[1]{\left\lceil #1 \right \rceil}
\newcommand{\floor}[1]{\left\lfloor #1 \right \rfloor}

\newcommand{\Hp}{\mathbb{H}^p(T)}
\newcommand{\Hpz}{\mathbb{H}^p_0(T)}

\renewcommand\root{\operatorname{o}}

\allowdisplaybreaks

\begin{document}
\date{\today}

\title[Shifts on $\Hp$] 
 {The Forward and Backward Shift on the Hardy Space of a Tree}

\author[A. \'Angeles-Romero]{Ad\'an \'Angeles-Romero}
\address{Universidad Aut\'onoma del Estado de Hidalgo\\ Pachuca, Mexico}
\email{adan.ar@hotmail.com }
\author[R.~A. Mart\'inez-Avenda\~no]{Rub\'en A. Mart\'inez-Avenda\~no}
\address{Departamento Acad\'emico de Matem\'aticas \\ Instituto
  Tecnol\'ogico Aut\'onomo de M\'exico \\ Mexico City, Mexico}
\email{ruben.martinez.avendano@gmail.com}


\thanks{The second author's research is partially supported by the
  Asociaci\'on Mexicana de Cultura A.C}

\begin{abstract}
  In this paper we initiate the study of the forward and backward shifts on the Hardy space of a tree and the little Hardy space of a tree. In particular, we investigate when these shifts are bounded, find the norm of the shifts if they are bounded, characterize the trees in which they are an isometry, compute the spectrum in some concrete examples, and completely determine when they are hypercyclic.
\end{abstract}

\maketitle

\section{Introduction}\label{sec_intro}

The study of Banach spaces of functions whose domains are trees, graphs, or discrete sets, has been studied for several decades now. The motivation for this study was, originally, to see what properties of harmonic or analytic functions on the unit disk, were preserved under a discretization of the domain (see, for example, \cite{Ca1,Ca2}). In particular, Cohen and Colonna \cite{CoCo1} studied the Bloch space of harmonic functions on a homogeneous tree and showed several interesting properties these functions have. In \cite{CoCo2}, they also showed how to embed certain homogeneous trees in the unit disk in such a way that the restriction of bounded harmonic functions on the disk are harmonic on the tree.

It is natural to endow the spaces of functions defined on a tree with norms that give them a Banach space structure. A natural norm on these spaces is, of course, the $L^p$ norm and this has been done in several places (see \cite{AgBeCoMaSl} for a list of some classical results). Perhaps the first example of a norm on the space of functions on a tree which was inspired by a classical space of holomorphic functions on the disk is the one introduced by Colonna and Easley in \cite{CoEa1}, where they study the so-called Lipschitz space of a tree. Furthermore, they investigate the properties of multiplication operators on this space, and on an importat subspace, the little Lipschitz space of a tree. (See also the paper \cite{AlCoEa}, where composition operators on the Lipschitz space of the tree are studied by Allen, Colonna and Easley.)

Muthukumar and Ponnusamy defined in \cite{MuPo1} the Hardy space of a tree, in analogy with the classical Hardy space of functions on the unit disk (see also \cite{KoMaTa} for an alternative viewpoint). They also introduce the little Hardy space of a tree, which does not have an analogous space in the classical case. In their paper, they also study properties of multiplication operators on the Hardy space of the tree and in \cite{MuPo2,MuPo3} they study the composition operators in these spaces.

On the other hand, the study of shift operators on different Banach spaces has produced many important results in operator theory and has been a source of examples and counterexamples to several questions in this area. Jab\l{}o\'{n}ski, Jung and Stochel, in the paper \cite{JaJuSt}, defined the forward shift on the $L^p$ space of functions defined on an infinite (rooted) tree and gave an explicit form for its adjoint, called the backward shift and studied further by the second author in \cite{Ma}. In the paper \cite{MaRi}, the second author and Rivera-Guasco studied properties of the forward and backward shifts defined on the Lipschitz and little Lipschitz space of a tree. We should point out that in \cite{AlJa}, Allen and Jackson study the {\em differentiation operator} on several discrete function spaces of a tree (including the Hardy space) and to define the differentiation operator they use a self-map of the tree which they call the {\em backward shift} and which, when used in a composition operator, would ``almost'' coincide with the operator we call the {\em forward shift} in this paper. We will not use their notation here.

In the present paper, we initiate the study of the forward and backward shifts on the Hardy space of a tree and on the little Hardy space of a tree. We investigate when these operators are bounded, when they are isometries, obtain the spectrum in some particular cases and characterize their hypercyclicity. We should point out that hypercyclicity of shift operators in spaces of functions on trees has been studied before, and in particular our results on hypercyclicity of the backward shift should be compared to those in \cite{GrPa,Ma,MaRi}. (See also the paper \cite{AlCoMaPo} for analogous results for composition operators on more general spaces of functions defined on discrete spaces.)

This paper is organized as follows. After giving in Section~\ref{sec_prelim} the definitions and notation that we will use throughout this paper, in Section~\ref{sec_forward} we completely characterize when the forward shift is bounded on the Hardy space and on the little Hardy space, and we give an explicit expression for its norm; see Subsection~\ref{subsec_boundS}. We also characterize those trees for which the forward shift is an isometry in Subsection~\ref{subsec_isometryS}. The most interesting question here is the question of what is the spectrum of the forward shift; after showing that for several cases the spectrum is the closed unit disk, we show an example of a tree where the spectrum is larger than the closed unit disk in Subsection~\ref{subsec_specS}. Lastly, in Subsection~\ref{subsec_hyperS} we show that the forward shift can never be hypercyclic. In Section~\ref{sec_backward}, we do the analogous results for the backward shift: we characterize when it is bounded, both on the Hardy space of the tree and on the little Hardy space of the tree and we obtain an expression for the norm of the operator in Subsection~\ref{subsec_boundedB}. We compute the spectrum for some particular examples, aided by some partial results on the spectrum of the backward shift in Subsection~\ref{subsec_specB}. We finish by showing in Subsection~\ref{subsec_hyperB} a complete characterization of the hyperciclicity of the backward shift on the little Hardy space of the tree.

\section{Preliminaries}\label{sec_prelim}

In this section we introduce the notation and definitions that will be used throughout the paper. As is customary, we will denote by $\N$, $\N_0$, $\R$, $\D$ and $\C$ the set of positive integers, the set of nonnegative integers, the set of real numbers, the set of complex numbers of modulus less than one, and the set of complex numbers, respectively.

\subsection{Trees, parents and children}
In this paper, we will consider infinite graphs which are locally finite; i.e., graphs in which the set of vertices is {\em countable} infinite but every vertex has a finite number of adjacent vertices. Furthermore, we consider only infinite trees; i.e., infinite graphs which are connected and which have no cycles. We will refer to infinite and locally finite trees simply as {\em trees}. 

We will always assume that each tree has a distinguished vertex, called the {\em root} and denoted by the symbol $\root$. We also use the usual distance for graphs: that is, given vertices $u$ and $v$ we denote by $\dist(u,v)$ the number of edges in the unique path between $u$ and $v$. For brevity's sake, we also define $|v|:=\dist(\root,v)$.  We will say that a vertex $v$ is in level $n$ if $|v|=n$. The number of vertices in level $n$ will be denoted by $\gamma(n)$.

Given a vertex $v \neq \root$, there is a unique vertex $w$ in the path between $v$ and the root $\root$ with $|w|=|v|-1$, we call this vertex the {\em parent} of $v$ and we denote it by $w=\parent{v}$. Inductively, for any integer $n \geq 2$,  we define the {\em $n$-parent} of $v$ as $\nparent{v}{n}:=\parent{\nparent{v}{n-1}}$, if $\nparent{v}{n-1} \neq \root$. 

We also denote by $\Chi{v}$, the {\em children} of $v$, the set of all vertices $u$ such that $v=\parent{u}$; observe that this set may be empty, in which case we say that $v$ is a leaf. If a tree has no leaves, we say the tree is leafless. Also, for each $n \in \N$ we define the set $\nChi{v}{n}$ as the set of all vertices $w$ such that $\nparent{w}{n}=v$; we refer to this set as the set of {\em $n$-children} of $v$.

For every vertex $w$, we define the {\em sector} $S_w$ as the set composed of $w$ and all of its $n$-children for every $n \in \N$. 

We finish this subsection with two important definitions that will be used throughout this paper.

\begin{itemize}
\item If $n \in \N$, we define $\gamma(n,v)$ to be the cardinality of $\nChi{v}{n}$. With this notation observe that $\gamma(1,v)$ is the number of children of the vertex $v$ and $\gamma(n,\root)=\gamma(n)$; i.e., the number of vertices at level $n$.

\item  If $n \in \N$ and $r \in \N_0$, we define $K(n,r):=\max\{\gamma(n,v) \, : \, |v|=r \}$. Observe that $K(n,r)$ is the number of $n$-children of the vertex that has the greatest number of $n$-children among the vertices at level $r$.
  \end{itemize}

 \subsection{Hardy space of the tree}

We are interested in complex-valued functions whose domain is the set of vertices of the tree. We denote the set of all such functions as $\calF$. For these functions, we make the following definition.

\begin{defi}
  Let $T$ be a tree and let $1 \leq p < \infty$. For $f \in \calF$ and for each $n \in \N_0$, we define
  \[
    M_p(n,f):=\left(\frac{1}{\gamma(n)} \sum_{|v|=n} |f(v)|^p\right)^{\frac{1}{p}}.
  \]
  That is, $M_p(n,f)$ is the ``$p$-th mean'' of the values of $|f(\cdot)|$ over all vertices at the level $n$.
 For economy of notation, we will write $M_p^p(n,f)$ instead of $(M_p(n,f))^p$.
\end{defi}

We now define the objects of study of this paper.

\begin{defi}
  Let $T$ be a tree and let $1 \leq p < \infty$  The {\em Hardy space of the tree}, denoted by $\Hp$, is the set of functions $f \in \calF$ such that there exists a positive $C \in \R$ with
  \[
    M_p(n,f) \leq C
  \]
  for all $n \in \N_0$.

  The  {\em little Hardy space of the tree}, denoted by $\Hpz$, is the set of functions $f \in \calF$ such that
   \[
     \lim_{n\to \infty }  M_p(n,f)=0.
   \]
 \end{defi}
 
 In \cite{MuPo1}, Muthukumar and Ponnusamy defined these spaces for regular trees (i.e., trees where every vertex has the same degree). They show that $\Hp$ is a Banach space when endowed with the norm
 \[
   \| f \| := \sup\left\{ M_p(n,f) \, : \, n \in \N_0 \right\},
 \]
 and that $\Hpz$ is a closed subspace of $\Hp$ and hence also a Banach space.

It is worth noting that Muthukumar and Ponnusamy's proofs that $\Hp$ and $\Hpz$ are Banach spaces apply not only for the case of nonregular trees which is what we will use in this paper, but for any graph with countably infinite vertices, which is locally finite and which has a distinguished vertex $\root$. We will not repeat their proofs here, but we leave it to the interested reader to check that indeed their proofs are valid in this more general context, essentially word for word.

Also, we should point out that none of the spaces $\Hp$ or $\Hpz$ are Hilbert spaces with the given norm, for any value of $p$, since it can be easily seen that this norm does not satisfy the parallelogram identity.

Lastly, it is easy to see that $\Hpz$ is a separable space. Indeed, observe that the set of functions of finite support,
\[
X:=\left\{ f \in \calF \, : \, f(v)\neq 0 \text{ for finitely many vertices } v \right\}
\]
is dense in $\Hpz$ and the set of functions in $X$ with (complex) rational values is countable and dense in $X$. On the other hand, $\Hp$ is not separable, as the set of functions in $\Hp$ with values in the set $\{ 0, 1\}$ and which are constant at each level is uncountable and the distance between any two (different) such functions is equal to $1$.

\subsection{Hypercyclicity} An operator $T$ on a Banach space $\calB$ is called {\em hypercyclic} if there exists $x \in \calB$ such that the orbit
\[
  \{ x, Tx, T^2x, T^3 x, \dots \}
  \]
is dense in $\calB$. Such a vector $x$ is also called a {\em hypercyclic vector}. Observe that if $T$ is hypercyclic, then the space $\calB$ is necessarily separable. It can be shown (see, for example, \cite{GrPe}) that if an operator is hypercyclic, the set of its hypercyclic vectors is a dense $G_\delta$-set in $\calB$. The main tool used to prove hypercyclicity in this paper is the Kitai-Gethner-Shapiro Criterion, which can be found, for example, in \cite[pp.74--75]{GrPe}

\section{The Forward Shift}\label{sec_forward}

In this section we study the forward shift on $\Hp$ and on $\Hpz$. We investigate when it is bounded, we completely characterize the trees where the forward shift is an isometry, find the spectrum in several cases, and show that it is never hypercyclic. First, we give the definition of the forward shift, originally given in \cite{JaJuSt}.

\begin{defi}
  Let $T$ be a tree. The {\em (forward) shift} $S: \calF \to \calF$ is defined as
  \[
    (Sf)(v)=\begin{cases}
      f(\parent{v}), & \text{ if } v \neq \root, \\
            0, & \text{ if } v =\root.
      \end{cases}
    \]
  \end{defi}

  It is clear that $S$ is a linear transformation on the vector space $\calF$.

\subsection{Boundedness of $S$}\label{subsec_boundS} The following theorem characterizes the boundedness of the forward shift on $\Hp$.

  \begin{theorem}\label{th:SboundedHp}
    Let $T$ be a tree and let $1 \leq p < \infty$. The forward shift $S$ is bounded on $\Hp$ if and only if the set
    \[
      \left\{ \frac{K(1,n-1) \gamma(n-1)}{\gamma(n)} \, : \, n \in \N \right\}
    \]
    is bounded. In this case,
    \[
      \| S \| = \left( \sup \left\{ \frac{K(1,n-1) \gamma(n-1)}{\gamma(n)} \, : \, n \in \N \right\} \right)^{\frac{1}{p}}.
    \]    
    \end{theorem}
    \begin{proof}
      Let us start by noticing that, for any $f \in \calF$ and any $n \in \N$ we have
      \begin{align*}
        M_p^p(n,Sf) & = \frac{1}{\gamma(n)} \sum_{|w|=n} |(Sf)(w)|^p \\
        & =  \frac{1}{\gamma(n)} \sum_{|w|=n} |f(\parent{w})|^p \\
        & =  \frac{1}{\gamma(n)} \sum_{|w|=n-1} \gamma(1,w) |f(w)|^p \\
        & \leq \frac{K(1,n-1)}{\gamma(n)} \sum_{|w|=n-1} |f(w)|^p \\
        & = \frac{K(1,n-1) \gamma(n-1)}{\gamma(n)} M_p^p(n-1,f),
      \end{align*}
      since $\gamma(1,w) \leq K(1,n-1)$ for every vertex $w$ with $|w|=n-1$. Hence
      \begin{equation}\label{eq:Sinequality}
        M_p^p(n,Sf) \leq \frac{K(1,n-1) \gamma(n-1)}{\gamma(n)} M_p^p(n-1,f).
        \end{equation}

        Now, let us assume that the set
        \[
          \left\{ \frac{K(1,n-1) \gamma(n-1)}{\gamma(n)} \, : \, n \in \N \right\}
        \]
        is bounded and let
        \[
          C:=\sup \left\{ \frac{K(1,n-1) \gamma(n-1)}{\gamma(n)} \, : \, n \in \N \right\}.
        \]
        From Equation \eqref{eq:Sinequality}, we have that
        \[
          M_p^p(n,Sf) \leq C M_p^p(n-1,f)
        \]
        for all $n \in \N$. Now, if $f \in \Hp$ it follows that
         \[
          M_p^p(n,Sf) \leq C \| f \|^p
        \]
        for all $n \in \N$ and also for $n=0$ since $M_p^p(0,Sf)=0$. This implies that $Sf \in \Hp$ and that
        \[
          \| Sf \|^p \leq C \| f \|^p.
        \]
        Hence $S$ is bounded. We also obtain in this case that \( \| S \| \leq C^{1/p} \).

        Let us assume now that $S$ is bounded. Choose a sequence of vertices $(w_n)$ such that $|w_n|=n$ and such that $\gamma(1,w_n)=K(1,n)$ for each $n \in \N_0$, that is, we choose a vertex at each level with the maximum number of children at that level. We define a function $f \in \calF$ as
        \[
          f(v):=\begin{cases}
            (\gamma(n))^{\frac{1}{p}}, & \text{ if } v=w_n \text{ for some } n \in \N_0,\\
            0, & \text{ in every other case}.
            \end{cases}
        \]
        For $n \in \N_0$, we have
        \[
          M_p^p(n,f)=\frac{1}{\gamma(n)} \sum_{|v|=n} |f(v)|^p = \frac{1}{\gamma(n)} |f(w_n)|^p =1,
        \]
        and hence $f \in \Hp$ and $\| f \|=1$. On the other hand, for every $n \in \N$ we have
        \begin{align*}
          M_p^p(n,Sf) &= \frac{1}{\gamma(n)} \sum_{|v|=n} | (Sf)(v)|^p \\
          &= \frac{1}{\gamma(n)} \sum_{|v|=n} | f(\parent{v})|^p \\
          &= \frac{1}{\gamma(n)} \sum_{v \in \Chi{w_{n-1}}} | f(\parent{v})|^p \\
          &= \frac{1}{\gamma(n)} \gamma(1, w_{n-1}) |f(w_{n-1})|^p \\
          &= \frac{1}{\gamma(n)} \gamma(1, w_{n-1}) |\gamma(n-1)^{\frac{1}{p}}|^p \\
          &= \frac{K(1,n-1) \gamma(n-1)}{\gamma(n)}.
        \end{align*}
        Since the operator $S$ is assumed to be bounded, it follows that $Sf \in \Hp$ and hence $M_p(n,Sf) \leq \|S f \|$ for all $n \in \N$. Therefore
        \[
          \frac{K(1,n-1) \gamma(n-1)}{\gamma(n)} \leq \|S f \|^p,
        \]
        for all $n \in \N$. It follows that the set
        \[
          \left\{ \frac{K(1,n-1) \gamma(n-1)}{\gamma(n)} \, : \, n \in \N \right\}
        \]
        is bounded and
        \[
          C=\sup\left\{ \frac{K(1,n-1) \gamma(n-1)}{\gamma(n)} \, : \, n \in \N \right\} \leq \|S f \|^p.
        \]
        Since $\| f \|=1$ we have that $C \leq \|S f \|^p \leq \|S\|^p$, from which we obtain $\| S \|= C^{1/p}$, as desired. This finishes the proof.
      \end{proof}

          The result above can be extended to $\Hpz$.

  \begin{theorem}\label{th:SboundedHpz}
    Let $T$ be a tree and let $1 \leq p < \infty$. The forward shift $S$ is bounded on $\Hpz$ if and only if the set
    \[
      \left\{ \frac{K(1,n-1) \gamma(n-1)}{\gamma(n)} \, : \, n \in \N \right\}
    \]
    is bounded. In this case,
    \[
      \| S \| = \left( \sup \left\{ \frac{K(1,n-1) \gamma(n-1)}{\gamma(n)} \, : \, n \in \N \right\} \right)^{\frac{1}{p}}.
    \]    
    \end{theorem}
    \begin{proof}
      First assume that
      \[
        \left\{ \frac{K(1,n-1) \gamma(n-1)}{\gamma(n)} \, : \, n \in \N \right\}
      \]
      is bounded. Then, by Theorem \ref{th:SboundedHp}, the operator $S$ is bounded on $\Hp$. Hence, to show that $S$ is bounded on $\Hpz$ it suffices to show that if $f \in \Hpz$ then $S f \in \Hpz$. As before, let
        \[
          C:=\sup \left\{ \frac{K(1,n-1) \gamma(n-1)}{\gamma(n)} \, : \, n \in \N \right\}.
        \]
        By Equation \eqref{eq:Sinequality}, we have that for all $n \in \N$
        \[
          M_p^p(n,Sf) \leq C M_p^p(n-1,f)
        \]
        and hence, since $M_p^p(n-1,f) \to 0$ as $n \to \infty$, it follows that $ M_p^p(n,Sf) \to 0$ as $n \to \infty$, showing that $Sf \in \Hpz$ and hence $S$ is bounded on $\Hpz$ with $\| S \| \leq C^{1/p}$.

        Now assume that $S$ is bounded on $\Hpz$. As we did in the proof of Theorem \ref{th:SboundedHp}, choose a sequence of vertices $(w_n)$ such that $|w_n|=n$ and such that $\gamma(1,w_n)=K(1,n)$ for each $n \in \N_0$. We now define a sequence of functions $(f_N) \subseteq \calF$ as follows. For each $N \in \N$, we set
        \[
          f_N(v):=\begin{cases}
            (\gamma(n))^{\frac{1}{p}}, & \text{ if } v=w_n \text{ for some } n \in \{ 1, 2, \dots, N\},\\
            0, & \text{ in every other case}.
            \end{cases}
          \]
          Observe that $M_p^p(n,f_N)=1$ if $0\leq n \leq N$ and $M_p^p(n,f_N)=0$ if $n \geq N$. Hence $f_N \in \Hpz$ and $\| f_N \|=1$ for all $N \in \N$. Also, as in the proof of Theorem \ref{th:SboundedHp}, observe that
          \[
            M_p^p(n, Sf_N)=\frac{K(1,n-1) \gamma(n-1)}{\gamma(n)},
          \]
          for $1 \leq n \leq N+1$, and $M_p^p(n, Sf_N)=0$ if $n > N+1$. Therefore, for all $n\in \N$ with $1 \leq n \leq N+1$ we have
          \begin{align*}
 \frac{K(1,n-1) \gamma(n-1)}{\gamma(n)} 
            &=  M_p^p(n,Sf_N) \\
            &\leq \sup \left\{ M_p^p(n,S f_N) \, : \, n \in \N_0  \right\} \\
            &= \| S f_N \|^p \\
            &\leq  \| S \|^p.
          \end{align*}
          Since $N$ is arbitrary, it follows that
          \[
            \frac{K(1,n-1) \gamma(n-1)}{\gamma(n)} \leq \| S \|^p
          \]
          for all $n \in \N$ and hence the set
          \[
      \left\{ \frac{K(1,n-1) \gamma(n-1)}{\gamma(n)} \, : \, n \in \N \right\}
    \]
    is bounded and
    \[
      C=\sup \left\{ \frac{K(1,n-1) \gamma(n-1)}{\gamma(n)} \, : \, n \in \N \right\} \leq \| S \|^p.
    \]
Since we already knew that $\|S \| \leq C^{1/p}$, this finishes the proof.
      \end{proof}

      Observe that for every tree, $\gamma(n) \leq \gamma(n-1) K(1,n-1)$ for all $n \in \N$ since, to obtain the vertices at the level $n$, we add to each vertex at level $n-1$ at most $K(1,n-1)$ vertices. Hence the norm of $S$ is always greater or equal to $1$. In fact, it clearly follows that $\|S\|=1$ if and only if $\gamma(n) = \gamma(n-1) K(1,n-1)$ for every $n \in \N$; that is, for each $n \in \N$, the number of children of all vertices at level $n$ are the same. We will say more about this case below (cf. Example~\ref{ex:Sisometry} and Proposition~\ref{prop:Sisometry}). We now show an example where the norm is strictly larger than $1$; in fact, we can find trees where the norm of $S$ is as large as we want.

      \begin{example}\label{ex:treek}
        Let $k \in \N$, with $k \geq 2$ and let $T$ be the tree defined as follows: start with the root $\root$ and add $k$ children to it. Choose one of these children, say $v_1$, and add $k$ children to it, and add only one child to each of the remaining $k-1$ vertices at level $1$. Choose one of the children of $v_1$, say $v_2$, add $k$ children to it, and add only one child to each of the $2k-2$ remaining vertices at level $2$. Choose one the children of $v_2$, say $v_3$, and add $k$ children to it, and add only one child to each of the remaining $3k-3$ vertices. Continue inductively. Then, the norm of $S$ in this tree is $k^{1/p}$.
        \end{example}
        \begin{proof}
          Observe that for all $n \in \N_0$ we have $\gamma(n)=nk-(n-1)$ and $K(1,n-1)=k$. Hence
          \[
            \sup \left\{ \frac{K(1,n-1) \gamma(n-1)}{\gamma(n)} \, : \, n \in \N \right\} = \sup \left\{ \frac{k \left( (n-1)k-(n-2)\right)}{nk-(n-1)} \, : \, n \in \N  \right\} = k
          \]
          from which the result follows.
          \end{proof}

          On the other hand, there are trees where $S$ is unbounded, as the following example shows.
          \begin{example}
            Construct a tree $T$ as follows. Let $\root$ be the root and add $2$ children to it. Pick one these two vertices, and add $3$ children to it, and just $1$ child to the other vertex. Inductively, at each level $n$, pick a vertex and add $n+2$ children to it, and add just $1$ child to each of the other vertices. Then $S$ is unbounded on $\Hp$ and on $\Hpz$.
          \end{example}
          \begin{proof}
            Observe that for all $n \in \N$ we have $\gamma(n)=\gamma(n-1)+n$. Since $\gamma(0)=1$ this implies that
            $\gamma(n)=\frac{n(n+1)}{2}+1$. Also note that $K(1,n)=n+2$. This implies that
            \[
              \frac{K(1,n-1) \gamma(n-1)}{\gamma(n)} = \frac{(n+1) \left(\frac{(n-1) n}{2}+1 \right)}{\frac{n(n+1)}{2}+1},
            \]
            which goes to infinity as $n$ goes to infinity. Hence, by Theorems~\ref{th:SboundedHp} and ~\ref{th:SboundedHpz} the operator $S$ is unbounded on $\Hp$ and on $\Hpz$.
          \end{proof}

      The following result will be useful when we compute the spectrum of $S$.
      \begin{theorem}\label{th:normSm}
        Let $T$ be a tree and let $1 \leq p < \infty$. If $S$ is bounded on $\Hp$ then
        \[
          \| S^m \| = \left( \sup \left\{ \frac{K(m,n-m) \gamma(n-m)}{\gamma(n)} \, : \, n \in \N, n \geq m \right\}\right)^{\frac{1}{p}}
        \]
        for each $m \in \N$.
        \end{theorem}
        \begin{proof}
          Let $m \in \N$ be fixed. Observe that for each $f \in \calF$ we have
          \[
            (S^m f)(v) = \begin{cases}
              f(\nparent{v}{m}), & \text{ if } v \in V^m,\\
                            0, & \text{ if } v \notin V^m,
              \end{cases}
            \]
            where $V^m$ is the set of vertices of $T$ which have an $m$-parent. Then, for $n \geq m$ we have
        \begin{align*}
        M_p^p(n,S^mf) & = \frac{1}{\gamma(n)} \sum_{|w|=n} |(S^mf)(w)|^p \\
        & =  \frac{1}{\gamma(n)} \sum_{|w|=n} |f(\nparent{w}{m})|^p \\
        & =  \frac{1}{\gamma(n)} \sum_{|w|=n-m} \gamma(m,w) |f(w)|^p \\
        & \leq \frac{K(m,n-m)}{\gamma(n)} \sum_{|w|=n-m} |f(w)|^p \\
        & = \frac{K(m,n-m) \gamma(n-m)}{\gamma(n)} M_p^p(n-m,f).
        \end{align*}
        Let
        \[
          C_m:= \sup \left\{ \frac{K(m,n-m) \gamma(n-m)}{\gamma(n)} \, : \, n \in \N, n \geq m \right\}
        \]
        Since $M_p^p(n,S^mf)=0$ if $0 \leq n < m$, it follows that
        \[
          M_p^p(n,S^mf) \leq C_m \| f \|^p,
        \]
        and hence
        \[
          \| S^m f \|^p \leq  C_m \| f \|^p,
        \]
        and therefore $\| S^m \| \leq C_m^{1/p}$.

        Now, let $( w_n)$ be a sequence of vertices such that $|w_n|=n$ and $\gamma(m,w_n)=K(m,n)$ for each $n \in \N_0$; that is, $w_n$ is the vertex at level $n$ with the largest number of $m$-children. Define $f \in \calF$ as
        \[
          f(v):=\begin{cases}
            (\gamma(n))^{\frac{1}{p}}, & \text{ if } v=w_n \text{ for some } n \in \N_0,\\
            0, & \text{ in every other case}.
            \end{cases}
        \]
        As before, $\| f \|=1$ and for $n \geq m$ we have
        \begin{align*}
          M_p^p(n,S^mf) &= \frac{1}{\gamma(n)} \sum_{|v|=n} | (S^mf)(v)|^p \\
          &= \frac{1}{\gamma(n)} \sum_{|v|=n} | f(\nparent{v}{m})|^p \\
          &= \frac{1}{\gamma(n)} \sum_{v \in \Chi{w_{n-m}}} | f(\nparent{v}{m})|^p \\
          &= \frac{1}{\gamma(n)} \gamma(m, w_{n-m}) |f(w_{n-m})|^p \\
          &= \frac{1}{\gamma(n)} \gamma(m, w_{n-m}) |\gamma(n-m)^{\frac{1}{p}}|^p \\
          &= \frac{K(m,n-m) \gamma(n-m)}{\gamma(n)}.
        \end{align*}        
        As before, since $S^m$ is bounded, it follows that
        \[
          \frac{K(1,n-m) \gamma(n-m)}{\gamma(n)} = M_p^p(n,S^mf) \leq \| S^m f \|^p \leq \| S^m \|^p
        \]
        and hence that $C_m \leq \| S^m \|^p$. Therefore
        \[
          \| S^m \| = C_m^{\frac{1}{p}},
        \]
        as desired.        
          \end{proof}

          \subsection{When is $S$ an isometry?}\label{subsec_isometryS} We start this investigation by looking at an example.  Observe that it is possible for the set of degrees of the tree to be unbounded while $S$ is an isometry.
      
      \begin{example}\label{ex:Sisometry}
        Let $(s_n)$ be a sequence of natural numbers. Let $T$ be the tree such that $\gamma(1,v)=s_{|v|+1}$ for all $v$. Then $S$ is an isometry on $\Hp$.
      \end{example}
      \begin{proof}
        First observe that $\gamma(n)=s_1 s_2 s_3 \cdots s_n$ for all $n \in \N$. Therefore, for all $n \in \N_0$ we have
        \begin{align*}
          M_p^p(n+1,Sf) &= \frac{1}{\gamma(n+1)} \sum_{|v|=n+1} |(Sf)(v)|^p \\
                        &= \frac{1}{\gamma(n+1)} \sum_{|v|=n+1} |f(\parent{v})|^p \\
                        &= \frac{1}{\gamma(n+1)} \sum_{|v|=n} s_{n+1}|f(v)|^p \\
                        &= \frac{1}{\gamma(n)} \sum_{|v|=n}|f(v)|^p \\
          &= M_p^p(n,f).
        \end{align*}
        Since $M_p^p(0,Sf)=0$, it follows that
        \[
          \sup \left\{ M_p^p (n, Sf) \, : \, n \in \N_0 \right\} =           \sup \left\{ M_p^p (n, f) \, : \, n \in \N_0 \right\},
        \]
        and hence $\|S f \|=\|f \|$, which finishes the proof.
      \end{proof}

      On the other hand, it turns out that if $S$ is an isometry, the tree has the above form.
      \begin{proposition}\label{prop:Sisometry}
        Let $T$ be a tree and assume that $S$ is an isometry on $\Hp$. Then there exists a sequence $(s_n)$ of natural numbers such that $\gamma(1,w)=s_{|w|+1}$ for all vertices $w$.
      \end{proposition}
      \begin{proof}
        Define the sequence $(s_n)$ as $s_n=\frac{\gamma(n+1)}{\gamma(n)}$. We will show that $(s_n)$ satisfies the desired conclusion.
        
        Let $w$ be a vertex and consider the function $f$ defined as the characteristic function of the vertex $w$. Let $n \in \N_0$. Observe that if $|w|=n$, then
        \[
          M_p^p(n,f)=\frac{1}{\gamma(n)} \sum_{|v|=n} |f(v)|^p = \frac{1}{\gamma(n)}, 
        \]
        but $M_p^p(n,f)=0$ if $|w|\neq n$. Hence $\| f \|^p= \frac{1}{\gamma(|w|)}$. But also observe that $Sf=\chi_{\Chi{w}}$, the characteristic function of the set $\Chi{w}$. Also, if $|w|=n-1$, then
        \[
          M_p^p(n,Sf)=\frac{1}{\gamma(n)} \sum_{|v|=n} |(Sf)(v)|^p = \frac{1}{\gamma(n)} \sum_{|v|=n} |\chi_{\Chi{w}}(v)|^p = \frac{1}{\gamma(n)} \gamma(1,w),
        \]
        and if $|w|\neq n-1$, then $M_p^p(n,Sf)=0$. Hence $\| S f\|^p= \frac{1}{\gamma(|w|+1)} \gamma(1,w)$. Since $S$ is an isometry, we must have
        \[
          \frac{1}{\gamma(|w|)}=\frac{1}{\gamma(|w|+1)} \gamma(1,w),
        \]
        from which we obtain that
        \[
          \gamma(1,w)=\frac{\gamma(|w|+1)}{\gamma(|w|)}=s_{|w|+1},
        \]
        and the result follows, since $w$ was arbitrary and $\gamma(1,w)$ is a natural number.
        \end{proof}

        Observe that the above two results also apply when $S$ is considered as an operator on $\Hpz$.
        
        \subsection{Spectrum of $S$}\label{subsec_specS} In this subsection we investigate the spectrum of the forward shift, both as an operator on $\Hpz$ and as an operator on $\Hpz$. We start with a simple observation (recall that $\sigma_p(S)$ denotes the set of eigenvalues).

        \begin{proposition}
          Let $T$ be a tree. If $T$ has a leaf, then $\sigma_p(S)=\{0\}$. If $T$ is leafless, then $\sigma_p(S)=\varnothing$.
        \end{proposition}
        \begin{proof}
First observe that if $\lambda \in \C$, with $\lambda \neq 0$, then $\lambda$ cannot be an eigenvalue of $S$. For suppose that $S f = \lambda f$: this implies that $0=(Sf)(\root)= \lambda f(\root)$, and hence $f(\root)=0$. We will show by induction that $f=0$. Suppose that $f(v)=0$ for every vertex $v$ with $|v|=n$, with $n \in \N_0$. If $w$ is a vertex with $|w|=n+1$, then $0 = f(\parent{w}) = (Sf)(w) = \lambda f(w)$ and hence $f(w)=0$, which shows $f$ is identically zero.

The only case left to check is $\lambda=0$. If $w$ is a leaf, then it is straightforward to check that $Sf=0$, where $f$ is the characteristic function of $w$. On the other hand, if $Sg=0$ for a nonzero function $g$, then $g(\parent{v})=0$ for each $v \neq \root$: this means that $g$ must be zero at every vertex which has a child. Since $g\neq 0$, it follows that $T$ must have a leaf. This finishes the proof.
          \end{proof}

          Observe that the result above is valid both in $\Hp$ and $\Hpz$, as is the following proposition.
          
          \begin{proposition}\label{pro:DsigmaS}
        Let $T$ be a tree. If $\lambda \in \D$, then $S-\lambda$ is not surjective. Hence $\conj{\D} \subseteq \sigma(S)$. 
      \end{proposition}
      \begin{proof}
        We will show that $\chi_{\root}$, the characteristic function of the root, is not in the range of $S-\lambda$ for $|\lambda|<1$. The result will follow from this.

        First of all, there is no function $f$ satisfying the equation $Sf=\chi_{\root}$, as can be seen by evaluating at the root. Now, assume $\lambda \neq 0$ and $(S-\lambda)f=\chi_{\root}$ for some function $f$. Evaluating at the root we obtain $0-\lambda f(\root)=1$ and hence $f(\root)=-\frac{1}{\lambda}$. Inductively, suppose that $f(v)=-\frac{1}{\lambda^{n+1}}$ for each vertex $v$ with $|v|=n$, with $n \in \N_0$. Let $w$ be a vertex with $|w|=n+1$. Evaluating $(S-\lambda)f=\chi_{\root}$ at $w$ we obtain $f(\parent{w})-\lambda f(w)=0$ and hence $-\frac{1}{\lambda^{n+1}} - \lambda f(w)=0$, from which it follows that $f(w)=-\frac{1}{\lambda^{n+2}}$, finishing the induction. But observe that then
        \[
          M_p (n, f)= \frac{1}{|\lambda^{n+1}|},
          \]
and therefore $f$ is not in $\Hp$ (nor in $\Hpz$) if $|\lambda|<1$. Hence $\chi_{\root}$ is not in the range of $S-\lambda$, as desired. The second part of the result follows, since $\sigma(S)$ is closed.
      \end{proof}
          
      A related result is the following.

      \begin{proposition} Let $T$ be a tree and $\lambda \in \C$ with $|\lambda|>1$. For each vertex $w$, the characteristic function $\chi_w$ is in the range of $S-\lambda$.
      \end{proposition}
      \begin{proof}
        Define the function $f_w \in \calF$ as follows:
        \[
          f_w(v)=\begin{cases}
            -\frac{1}{\lambda^{|v|-|w|+1}}, & \text{ if } v \in S_w\\
            0, & \text{ in any other case}.
  \end{cases}
\]
First observe that $((S-\lambda)f_w)(w)=f_w(\parent{w})-\lambda f_w(w) =0-\lambda (-\frac{1}{\lambda})=1=\chi_w(w)$. If $v \in S_w$ but $v \neq w$, then
\begin{align*}
  ((S-\lambda) f_w)(w) &= f_w(\parent{w})-\lambda f_w(w) \\ &=-\frac{1}{\lambda^{|\parent{v}|-|w|+1}} - (-\lambda) \frac{1}{\lambda^{|v|-|w|+1}}\\
                                                      &= -\frac{1}{\lambda^{(|v|-1)-|w|+1}} + \frac{1}{\lambda^{|v|-|w|}} \\
                                                      &= 0 \\
  &= \chi_w(v),
\end{align*}
while if $v \notin S_w$ we have $((S-\lambda)f_w)(v)=f_w(\parent{v})-\lambda f_w(v) = 0 - \lambda \cdot 0 = 0 =\chi_w(v)$. Hence $(S-\lambda)f_w=\chi_w$. The proof will be done by observing that for every $n > |w|$ we have
\[
  M_p^p(n,f_w)=\frac{1}{\gamma(n)} \gamma(n-|w|,w) \left( \frac{1}{|\lambda|^{n-|w|+1}} \right)^p \leq \left( \frac{1}{|\lambda|^{n-|w|+1}} \right)^p.
\]
But this goes to $0$ as $n \to \infty$ and hence $f \in \Hp$ and $f \in \Hpz$.
        \end{proof}

        In the case of $\Hpz$, since the set of functions of finite support is dense in $\Hpz$, the above result gives the following result (recall that $\sigma_{\text{ap}}(S)$ denotes the approximate point spectrum of $S$).
        \begin{corollary}
          Let $T$ be a tree and let $\lambda \in \C$, with $|\lambda|>1$. Then the range of $S-\lambda$ is dense in $\Hpz$. In particular, $\sigma(S) \cap (\C \setminus \conj{\D}) \subseteq \sigma_{\text{ap}} (S)$ in $\Hpz$.
        \end{corollary}

        The following examples show that, depending on the tree, the spectrum of $S$ may be the unit disk or may be larger than the unit disk.

        \begin{example}
          Let $T$ be a tree and assume $\|S \|=1$. Then $\sigma(S)=\conj{\D}$.
        \end{example}
        \begin{proof}
          This follows from $\|S\|=1$ and Proposition \ref{pro:DsigmaS}.
        \end{proof}
Of course, the hypothesis in the above example only hold when $S$ is an isometry (see Example \ref{ex:Sisometry} and Proposition \ref{prop:Sisometry}).

The following example shows that the spectrum may be the unit disk even if the norm is larger than $1$.

        \begin{example}
          Let $T$ be the tree in Example \ref{ex:treek}. Then $\sigma(S)=\conj{\D}$ and $\| S \|=k^{1/p}$.
        \end{example}
        \begin{proof}
          We showed in Example \ref{ex:treek} that  $\| S \|=k^{1/p}$. For the rest of the proof, we will first compute the spectral radius of $S$. Observe that $\gamma(n)=nk-(n-1)$ for all $n \in \N$ and $K(m,t)=mk-(m-1)$ for all $m \in \N$ and all $t \in \N_0$. For each fixed $m$ we have
          \begin{align*}
            \sup \left\{ \frac{K(m,n-m) \gamma(n-m)}{\gamma(n)} \, : \, n \geq m \right\} &= \sup \left\{\frac{(mk-(m-1)) ((n-m)k-(n-m-1)}{nk -(n-1)} \, : \, n \geq m \right\} \\
            &= (mk-(m-1) ) \sup \left\{\frac{((n-m)k-(n-m-1)}{nk -(n-1)} \, : \, n \geq m \right\} \\
            &=mk-(m-1).
\end{align*}
But observe that
\[
  \lim_{m \to \infty} (mk -(m-1))^{1/m}=1,
\]
and hence $\lim_{m\to\infty} \| S^m\|^{1/m} =1$. Therefore $\sigma(S) \subseteq \conj{\D}$. But now Proposition \ref{pro:DsigmaS} completes the claim.
          \end{proof}

          For all the examples we have seen so far, the spectrum of the shift equals to closed unit disk. Nevertheless, there are examples of trees where the spectral radius is larger than $1$.

\begin{example}
  Let $T$ be the tree defined as follows (see the picture below): start with the root $\root$ and add two children to it, say $u^1_1$ and $u^1_2$ (picture them ordered from left to right). We choose $\ceiling{\frac{2}{2}}=1$ of these vertices, the vertex $u_1^1$ and we add two children to it, say $u^2_1$ and $u^2_2$ (again, picture them ordered from left to right) and we add only one child to $u^1_2$, say $v^2_1$. Now, of the $3$ vertices at level $2$, we choose $\ceiling{\frac{3}{2}}=2$ of them, $u^2_1$ and $u^2_2$ and we add two children to each: from left to right, we call the children of $u^2_1$ by the names $u^3_1$, $u^3_2$ and the children of $u^2_2$ by the names $u^3_3$, $u^3_4$. We add only one child to the vertex $v^2_1$, and we call it $v^3_1$.

We now have $5$ vertices at level $3$: choose  $\ceiling{\frac{5}{2}}=3$ of them, the vertices $u^3_1$, $u^3_2$ and $u^3_3$ and add two children to each, say $u^4_1$,  $u^4_2$, $u^4_3$, $u^4_4$, $u^4_5$ and $u^4_6$; also, we add only one child to each of the vertices  $u^3_4$ and $v^3_1$, say $v^4_1$ and $v^4_2$, respectively.

Inductively, if we have $a_n$ vertices at level $n$, choose $\ceiling{\frac{a_n}{2}}$ of them, namely $u^n_j$ for $j=1, 2, \dots, \ceiling{\frac{a_n}{2}}$ and we add two children to each to obtain, from left to right, the vertices $u^{n+1}_j$, with $j=1, 2, \dots, 2 \ceiling{\frac{a_n}{2}}$. We only add one child to the remaining vertices $u^n_j$, for $j \geq \ceiling{\frac{a_n}{2}}+1$ and $v^n_j$, for all $j$, from left to right denoted by $v^{n+1}_j$, for $j=1, 2, \dots a_n - \ceiling{\frac{a_n}{2}}$.

We describe this process in the following picture. The dashed blue lines indicate the level (they are not edges), while the dotted black lines indicate that the process continues.

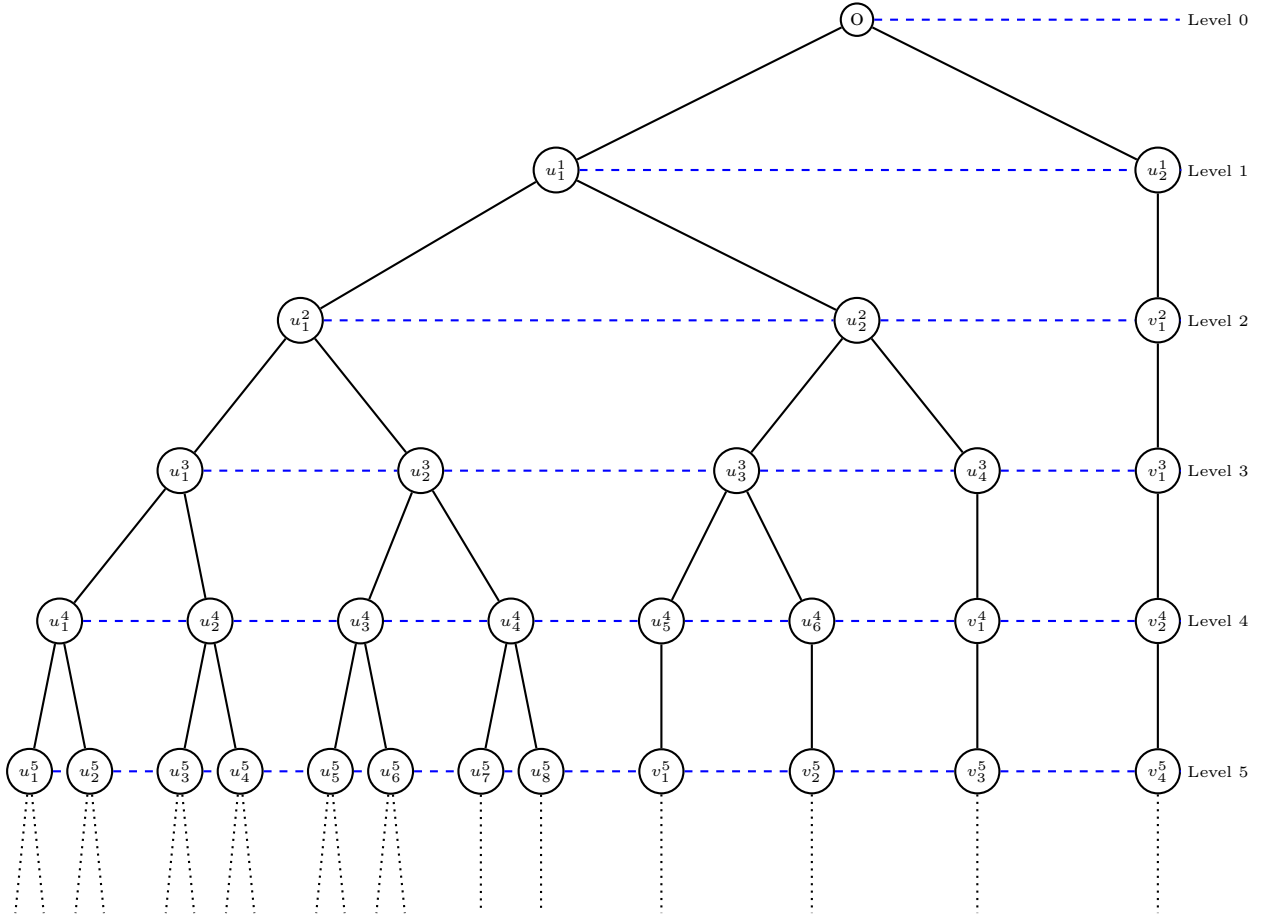
\begin{center}
	\begin{tikzpicture}[scale=2,thick]
		\tikzstyle{every node}=[minimum width=0pt, inner sep=2pt, circle]
			\draw (0.5,6) node[draw] (0) {$\root$};
			\draw (-1.5,5) node[draw] (1) { \tiny $u_1^1$};
			\draw (2.5,5) node[draw] (2) { \tiny $u_2^1$};
			\draw (-3.2,4) node[draw] (3) { \tiny $u^2_1$};
			\draw (0.5,4) node[draw] (4) { \tiny $u^2_2$};
			\draw (2.5,4) node[draw] (5) { \tiny $v^2_1$};
			\draw (-4,3) node[draw] (6) { \tiny $u^3_1$};
			\draw (-2.4,3) node[draw] (7) { \tiny $u^3_2$};
			\draw (-0.3,3) node[draw] (8) { \tiny $u^3_3$};
			\draw (1.3,3) node[draw] (9) { \tiny $u^3_4$};
			\draw (2.5,3) node[draw] (10) { \tiny $v^3_1$};
			\draw (-4.8,2) node[draw] (11) { \tiny  $u^4_1$};
			\draw (-3.8,2) node[draw] (12) { \tiny  $u^4_2$};
			\draw (-2.8,2) node[draw] (13) { \tiny  $u^4_3$};
			\draw (-1.8,2) node[draw] (14) { \tiny  $u^4_4$};
			\draw (-0.8,2) node[draw] (15) { \tiny  $u^4_5$};
			\draw (0.2,2) node[draw] (16) { \tiny  $u^4_6$};
			\draw (1.3,2) node[draw] (17) { \tiny  $v^4_1$};
			\draw (2.5,2) node[draw] (18) { \tiny  $v^4_2$};
                        \draw (-5,1) node[draw] (19) { \tiny  $u^5_1$};
			\draw (-4.6,1) node[draw] (20) { \tiny  $u^5_2$};
			\draw (-4,1) node[draw] (21) { \tiny  $u^5_3$};
			\draw (-3.6,1) node[draw] (22) { \tiny  $u^5_4$};
                        \draw (-3,1) node[draw] (23) { \tiny  $u^5_5$};
			\draw (-2.6,1) node[draw] (24) { \tiny  $u^5_6$};
			\draw (-2,1) node[draw] (25) { \tiny  $u^5_7$};
			\draw (-1.6,1) node[draw] (26) { \tiny  $u^5_8$};
                        \draw (-0.8,1) node[draw] (27) { \tiny  $v^5_1$};
			\draw (0.2,1) node[draw] (28) { \tiny  $v^5_2$};
			\draw (1.3,1) node[draw] (29) { \tiny  $v^5_3$};
			\draw (2.5,1) node[draw] (30) { \tiny  $v^5_4$};
			\draw (-5.1,0) node (31) {};
                        \draw (-4.9,0) node (32) {};
			\draw (-4.7,0) node (33) {};
                        \draw (-4.5,0) node (34) {};
			\draw (-4.1,0) node (35) {};
                        \draw (-3.9,0) node (36) {};
			\draw (-3.7,0) node (37) {};
                        \draw (-3.5,0) node (38) {};
			\draw (-3.1,0) node (39) {};
                        \draw (-2.9,0) node (40) {};
			\draw (-2.7,0) node (41) {};
                        \draw (-2.5,0) node (42) {};
			\draw (-2,0) node (43) {};
                        \draw (-1.6,0) node (44) {};
			\draw (-0.8,0) node (45) {};
                        \draw (0.2,0) node (46) {};
                        \draw (1.3,0) node (47) {};
                        \draw (2.5,0) node (48) {};
                        \draw (2.9,6) node (49) {\tiny Level 0};
                        \draw (2.9,5) node (50) {\tiny Level 1};
                        \draw (2.9,4) node (51) {\tiny Level 2};
                        \draw (2.9,3) node (52) {\tiny Level 3};
                        \draw (2.9,2) node (53) {\tiny Level 4};
                        \draw (2.9,1) node (54) {\tiny Level 5};
                        \draw  (0) edge[blue,dashed] (49);
                        \draw  (1) edge[blue,dashed] (2);
                        \draw  (2) edge[blue,dashed] (50);
                        \draw  (3) edge[blue,dashed] (4);
                        \draw  (4) edge[blue,dashed] (5);
                        \draw  (5) edge[blue,dashed] (51);
                        \draw  (6) edge[blue,dashed] (7);
                        \draw  (7) edge[blue,dashed] (8);
                        \draw  (8) edge[blue,dashed] (9);
                        \draw  (9) edge[blue,dashed] (10);
                        \draw  (10) edge[blue,dashed] (52);
                        \draw  (11) edge[blue,dashed] (12);
                        \draw  (12) edge[blue,dashed] (13);
                        \draw  (13) edge[blue,dashed] (14);
                        \draw  (14) edge[blue,dashed] (15);
                        \draw  (15) edge[blue,dashed] (16);
                        \draw  (16) edge[blue,dashed] (17);
                        \draw  (17) edge[blue,dashed] (18);
                        \draw  (18) edge[blue,dashed] (53);
                        \draw  (19) edge[blue,dashed] (20);
                        \draw  (20) edge[blue,dashed] (21);                        
                        \draw  (21) edge[blue,dashed] (22);                        
                        \draw  (22) edge[blue,dashed] (23);                        
                        \draw  (23) edge[blue,dashed] (24);                        
                        \draw  (24) edge[blue,dashed] (25);                        
                        \draw  (25) edge[blue,dashed] (26);                        
                        \draw  (26) edge[blue,dashed] (27);                        
                        \draw  (27) edge[blue,dashed] (28);                        
                        \draw  (28) edge[blue,dashed] (29);                        
                        \draw  (29) edge[blue,dashed] (30);                        
                        \draw  (30) edge[blue,dashed] (54);                        
                        \draw  (0) edge (1);
			\draw  (0) edge (2);
			\draw  (1) edge (3);
			\draw  (1) edge (4);
			\draw  (2) edge (5);
			\draw  (3) edge (6);
			\draw  (3) edge (7);
			\draw  (4) edge (8);
			\draw  (4) edge (9);
			\draw  (5) edge (10);
			\draw  (6) edge (11);
			\draw  (6) edge (12);
			\draw  (7) edge (13);
			\draw  (7) edge (14);
			\draw  (8) edge (15);
			\draw  (8) edge (16);
			\draw  (9) edge (17);
			\draw  (10) edge (18);
                        \draw  (11) edge (19);
                        \draw  (11) edge (20);
                        \draw  (12) edge (21);
                        \draw  (12) edge (22);
                        \draw  (13) edge (23);
                        \draw  (13) edge (24);
                        \draw  (14) edge (25);
                        \draw  (14) edge (26);
                        \draw  (15) edge (27);
                        \draw  (16) edge (28);
                        \draw  (17) edge (29);
                        \draw  (18) edge (30);
                        \draw  (19) edge[dotted] (31);
                        \draw  (19) edge[dotted] (32);
                        \draw  (20) edge[dotted] (33);
                        \draw  (20) edge[dotted] (34);
                        \draw  (21) edge[dotted] (35);
                        \draw  (21) edge[dotted] (36);
                        \draw  (22) edge[dotted] (37);
                        \draw  (22) edge[dotted] (38);
                        \draw  (23) edge[dotted] (39);
                        \draw  (23) edge[dotted] (40);
                        \draw  (24) edge[dotted] (41);
                        \draw  (24) edge[dotted] (42);
                        \draw  (25) edge[dotted] (43);
                        \draw  (26) edge[dotted] (44);
                        \draw  (27) edge[dotted] (45);
                        \draw  (28) edge[dotted] (46);
                        \draw  (29) edge[dotted] (47);
                        \draw  (30) edge[dotted] (48); 
                      \end{tikzpicture}
                    \end{center}
Then, the spectrum of $S$ on $\Hp$ is not the closed unit disk. In fact, the spectral radius of $S$ is $\frac{4}{3}$.  
 \end{example}
                \begin{proof}
            Observe that the number of children $a_n:=\gamma(n)$ at level $n \in \N_0$ obeys the recursive relation
            \[
              a_{n+1}= 2 \ceiling{\frac{a_n} {2}} + \floor{ \frac{a_n}{2}}             \]
            with $a_0=1$. This relation can be simplified to obtain
            \[
              a_{n+1}= \ceiling{\frac{3 a_n}{2}}.
            \]
            It can be seen (see \cite{OdWi}) that $a_n=\mathrm{O}\left((\frac{3}{2})^n\right)$ and, in fact,
            \[
              a_n=\floor{c \left( \frac{3}{2} \right)^n},
            \]
            where $c$ is a constant (it is known that $c\approx 1.62$).
            
Now, we claim that for each level $m \in \N$, the subgraph spanned by $u^m_1$ and $\ds\bigcup_{k=1}^{m+1}\nChi{u^m_1}{k}$ is a binary tree of depth $m+1$. We prove this by induction on $m$.

            If $m=1$, clearly the subgraph spanned by $u^1_1$, the set $\Chi{u^1_1}=\{u^2_1,u^2_2\}$ and the set $\nChi{u^1_1}{2}=\{u^3_1,u^3_2,u^3_3,u^3_4\}$ is a binary tree of depth $2$.

            Now, assume that we have shown that for some level $m$, the subgraph spanned by
            \[
               \{ u^m_1\} \cup  \bigcup_{k=1}^{m+1}\nChi{u^m_1}{k}
             \]
             is a binary tree of depth $m+1$. Observe that this set contains the set $\ds\{ u^{m+1}_1\} \cup \bigcup_{k=1}^{m}\nChi{u^{m+1}_1}{k}$, which by the induction hypothesis spans a binary tree of depth $m$. This implies that the set $\nChi{u^{m+1}_1}{m}$ contains exactly $2^m$ vertices, but the induction hypothesis also implies that the set $\nChi{u^{m+1}_2}{m}$ also contains $2^m$ vertices: it follows that the set $\nChi{u^{m+1}_1}{m}$  contains less than half of all vertices at the level $2m+1$. Hence, the construction of the tree implies that every vertex in the set $\nChi{u^{m+1}_1}{m}$ has two children, and therefore the set
            \[
              \{ u^{m+1}_1\} \cup \bigcup_{k=1}^{m+1}\nChi{u^{m+1}_1}{k}
            \]
            spans a binary tree of depth $m+1$. Furthermore, observe that the set $\nChi{u^{m+1}_1}{m+1}$ has $2^{m+1}$ vertices, and this set is inside the level $2m+2$, which has $a_{2m+2}$ vertices. It can easily be proved that $2^{m+1} \leq \frac{1}{2} a_{2m+2}$, therefore every vertex in the set $\nChi{u^{m+1}_1}{m+1}$ will have two children, by construction of the tree. Hence the set 
            \[
              \{ u^{m+1}_1\} \cup \bigcup_{k=1}^{m+2}\nChi{u^{m+1}_1}{k}
            \]
            spans a binary tree of depth $m+2$, completing the induction.

            The previous claim says that at each level $m$, there is a vertex who is the root of a binary subtree of depth $m+1$. Since the tree $T$ has the property that every vertex has either one or two children, it follows that $K(r,s)$, the greatest number of $r$-children of the vertices at level $s$ is less than $2^r$ and, in fact, $K(r,s)=2^r$ if $1 \leq r \leq s+1$.

            Now, fix $m \in \N$. Theorem~\ref{th:normSm} implies that
            \[
              \| S^m \| = \sup\left\{\frac{K(m,n-m) a_{n-m} }{a_n} \, : \, n \in \N, n \geq m \right\}.
            \]
Also, observe that
                \begin{align*}
                  \| S^m \| &= \sup\left\{\frac{K(m,n-m) a_{n-m} }{a_n} \, : \, n \in \N, n \geq m \right\} \\
                            &=\max\left\{
                              \sup\left\{\frac{K(m,n-m) a_{n-m} }{a_n} \, : \, m \leq n <2 m-1 \right\},  \sup\left\{\frac{K(m,n-m) a_{n-m} }{a_n} \, : \, n \geq 2m-1 \right\}
                              \right\}.
                \end{align*}
It is straightfoward to verify that $\frac{a_{n-m}}{a_n} \leq \left(\frac{2}{3}\right)^m$ and that, in fact, $\lim_{n \to \infty}\frac{a_{n-m}}{a_n} = \left(\frac{2}{3}\right)^m$. Since $K(m,n-m)= 2^{m}$ if $1 \leq m \leq n-m+1$ and $K(m,n-m) \leq 2^{m}$ otherwise, we have that if $m\leq n < 2m-1$, then
\[
  \frac{K(m,n-m) a_{n-m} }{a_n} \leq 2^m \left(\frac{2}{3}\right)^m = \frac{4^m}{3^m},
\]
and therefore
\[
   \sup\left\{\frac{K(m,n-m) a_{n-m} }{a_n} \, : \, n <2 m-1 \right\} \leq  \frac{4^m}{3^m}.
\]                          
On the other hand, since we have
\[
   \sup\left\{\frac{K(m,n-m) a_{n-m} }{a_n} \, : \, n \geq 2m-1 \right\} = 2^m  \sup\left\{\frac{a_{n-m} }{a_n} \, : \, n \geq 2m-1 \right\},
 \]
$\frac{a_{n-m}}{a_n} \leq \left(\frac{2}{3}\right)^m$ and  $\lim_{n \to \infty}\frac{a_{n-m}}{a_n} = \left(\frac{2}{3}\right)^m$, it follows that
     \[
   \sup\left\{\frac{K(m,n-m) a_{n-m} }{a_n} \, : \, n \geq 2m-1 \right\} = 2^m  \frac{2^m}{3^m}.
 \]
 Hence,
 \[
   \| S^m \| = \frac{4^m}{3^m}. 
 \]
 But by the spectral radius formula, the spectral radius of $S$ equals
 \[
   \lim_{m\to \infty} \| S^m \|^{\frac{1}{m}}=\frac{4}{3},
 \]
 and therefore the spectrum of $S$ is not the closed unit disk.
\end{proof}

What is the spectrum of the operator in the example above? We leave that question open for future research.

\subsection{Hypercyclicity of $S$}\label{subsec_hyperS} We finish this section by studying the hypercyclicity of the operator $S$. First observe that since $\Hp$ is not separable, $S$ cannot be hypercyclic on $\Hp$. It turns out that even though $\Hpz$ is separable, $S$ is not hypercyclic on $\Hpz$.

\begin{proposition}
Let $T$ be a tree and assume $S$ is bounded on $\Hpz$. Then $S$ is not hypercyclic on $\Hpz$.
\end{proposition}
\begin{proof}
  Assume, by contradiction, that $S$ is hypercyclic and let $f \in \Hpz$ be a hypercyclic vector. Let $g \in \Hpz$ be such that $g(\root)\neq 0$. Then, there exists $N \in \N$ such that $\| S^N f - g \| < |g(\root)|$. But since
  \[
    |g(\root)| = |(S^N f-g)(\root)| = M_p(0,S^N f - g) \leq \| S^N f - g \|,
  \]
  this is a contradiction.
\end{proof}

\section{The Backward Shift}\label{sec_backward}

In this section, we investigate the boundedness of the backward shift operator, its spectrum in several cases and its hypercyclicity on $\Hpz$. We start with the definition of the backward shift.

\begin{defi}
  Let $T$ be a tree. The {\em backward shift} $B:\calF \to \calF$ is defined as
  \[
    (Bf)(v)=\sum_{w \in \Chi{v}} f(w),
  \]
  where, as usual, if $\Chi{v}$ is empty, the sum is defined to be equal to $0$.
\end{defi}

Clearly $B$ is a linear transformation on $\calF$.

\subsection{Boundedness of $B$}\label{subsec_boundedB} The following theorem characterizes the boundedness of the backward shift on $\Hp$.

\begin{theorem}\label{th:BboundedHp}
  Let $T$ be a tree and let $1 \leq p < \infty$. The backward shift $B$ is bounded on $\Hp$ if and only if the set
  \[
    \left\{ \frac{(K(1,n))^{p-1} \gamma(n+1)}{\gamma(n)} \, : \, n \in \N_0 \right\}
  \]
  is bounded. In this case,
  \[
    \| B \| = \left(\sup \left\{ \frac{(K(1,n))^{p-1} \gamma(n+1)}{\gamma(n)} \, : \, n \in \N_0 \right\} \right)^{1/p}.
  \]
  \end{theorem}
  \begin{proof}
    First observe that for any $f \in \calF$ and $n \in \N_0$ we have
\begin{align*}
  M_p^p(n, Bf) &=\frac{1}{\gamma(n)} \sum_{|w|=n}|(Bf)(w)|^p \\
               &=\frac{1}{\gamma(n)} \sum_{|w|=n}\left|\sum_{v\in\Chi{w}} f(v) \right|^p\\
               &\leq \frac{1}{\gamma(n)} \sum_{|w|=n}\left(\sum_{v\in\Chi{w}} |f(v)| \right)^p\\
               &\leq \frac{1}{\gamma(n)}\sum_{|w|=n}\gamma(1,w)^{p-1}\sum_{v\in \Chi{w}}|f(v)|^p,
\end{align*}
where the last inequality follows from Jensen's inequality. Since $\gamma(1,w) \leq K(1,n)$ for every vertex $w$ with $|w|=n$, we have  
\begin{equation}\label{eq:Binequality}
  M_p^p(n, Bf)  \leq \frac{K(1,n)^{p-1}}{\gamma(n)} \sum_{|v|=n+1}|f(v)|^p = \frac{(K(1,n))^{p-1}\gamma(n+1)}{\gamma(n)}M_{p}^{p}(n+1,f).
\end{equation}
Assume now that the set
  \[
    \left\{ \frac{(K(1,n))^{p-1} \gamma(n+1)}{\gamma(n)} \, : \, n \in \N_0 \right\}
  \]
  is bounded and let
  \[
    C:= \sup \left\{ \frac{(K(1,n))^{p-1} \gamma(n+1)}{\gamma(n)} \, : \, n \in \N_0 \right\}.
  \]
  By equation \eqref{eq:Binequality} we obtain
  \[
    M_p^p(n,Bf) \leq C \ M_p^p(n+1,f),
  \]
  for all $n \in \N_0$. Now, if $f \in \Hp$ we have
  \[
    M_p^p(n,Bf) \leq C \| f \|^p,
  \]
  and hence $Bf \in \Hp$ and
  \[
    \| Bf \|^p \leq C \| f \|^p.
  \]
  Therefore $B$ is bounded. We also obtain in this case that $\| B \| \leq C^{1/p}$.

  Now, let us assume $B$ is bounded on $\Hp$. For every $n \in \N_0$, choose a vertex $w_n$ such that $|w_n|=n$ and $\gamma(1,w_n)=K(1,n)$. We define the function
  \[
    f(v)=\begin{cases}
      1,& \text{ if } v=\root,\\
      \left(\frac{\gamma(n+1)}{K(1,n)}\right)^{\frac{1}{p}}, & \text{ if } v \in \Chi{w_{n}} \text{ for some } n \in \N_0, \text{ and} \\
      0, & \text{ in any other case.}
      \end{cases}
    \]
    Observe thet $M_{p}^{p}(0,f)=1$ and that for all $n \in \N$ we have
\begin{align*}
  M_{p}^{p}(n,f)
  &=\frac{1}{\gamma(n)}\sum_{|v|=n}|f(v)|^p \\
  &= \frac{1}{\gamma(n)}\sum_{v\in \Chi{w_{n-1}}}|f(v)|^p \\
  &= \frac{K(1,n-1)}{\gamma(n)}\left|\left( \frac{\gamma(n)}{K(1,n-1)}\right)^{\frac{1}{p}}\right|^p \\
  &= 1.
\end{align*}
Hence  $f\in\Hp$ and $\| f \|=1$.

Observe that if $|v|=n$ but $v\neq w_n$, then
\[
  (Bf)(v)=\sum_{w \in \Chi{v}} f(w)=0.
\]

Hence, for every $n \in \N_0$, we have
\begin{align*}
  M_{p}^{p}(n,Bf)
  &=\frac{1}{\gamma(n)}\sum_{|v|=n}|(Bf)(v)|^{p} \\
  &=\frac{1}{\gamma(n)}|Bf(w_{n})|^{p}\\
  &= \frac{1}{\gamma(n)}\left|\sum_{u\in \Chi{w_{n}}}f(u)\right|^{p}\\
  &=\frac{1}{\gamma(n)}\left|K(1,n)\left(\frac{\gamma(n+1)}
  {K(1,n)}\right)^{\frac{1}{p}}\right|^{p}\\
  &=\frac{(K(1,n))^{p-1}\gamma(n+1)}{\gamma(n)}.
\end{align*}
Since $B$ is bounded, this implies that
\[
\frac{(K(1,n))^{p-1}\gamma(n+1)}{\gamma(n)}=  M_{p}^{p}(n,Bf) \leq \| B \|^p,
\]
and hence the set
  \[
    \left\{ \frac{(K(1,n))^{p-1} \gamma(n+1)}{\gamma(n)} \, : \, n \in \N_0 \right\}
 \]
  is bounded and
 \[
C=\sup \left\{ \frac{(K(1,n))^{p-1} \gamma(n+1)}{\gamma(n)} \, : \, n \in \N_0 \right\}  \leq \|B \|^p.
  \] 
This concludes the proof.
\end{proof}

The result can be extended to $\Hpz$.

\begin{theorem}\label{th:BboundedHpz}
  Let $T$ be a tree and let $1 \leq p < \infty$. The backward shift $B$ is bounded on $\Hpz$ if and only if the set
  \[
    \left\{ \frac{(K(1,n))^{p-1} \gamma(n+1)}{\gamma(n)} \, : \, n \in \N_0 \right\}
  \]
  is bounded. In this case,
  \[
    \| B \| = \left(\sup \left\{ \frac{(K(1,n))^{p-1} \gamma(n+1)}{\gamma(n)} \, : \, n \in \N_0 \right\} \right)^{1/p}.
  \]
  \end{theorem}
  \begin{proof}
   First assume the set
  \[
    \left\{ \frac{(K(1,n))^{p-1} \gamma(n+1)}{\gamma(n)} \, : \, n \in \N_0 \right\}
  \]
  is bounded. By Theorem~\ref{th:BboundedHp}, is is enough to show that if $f \in \Hpz$, then $Bf \in \Hpz$. By inequality \eqref{eq:Binequality}, we have
  \[
    M_p^p(n,Bf)\leq C \ M_p^p(n+1,f)
  \]
  for all $n \in \N_0$. Hence, if $f \in\Hpz$, then $Bf \in \Hpz$, as desired.

  For the other direction, assume that $B$ is bounded on $\Hpz$. For each $n \in \N_0$, we choose a vertex $w_n$ with $|w_n|=n$ and $\gamma(1,w_n)=K(1,n)$. We now define a sequence of functions $(f_N) \subseteq \calF$ as follows.
For each $N \in \N$ we set
  \[
    f_N(v)= \begin{cases}
      1, & \text{ if } v=\root,\\
      \left(\frac{\gamma(n+1)}{K(1,n)}\right)^{\frac{1}{p}}, & \text{ if } v \in \Chi{w_n} \text{ and } 0 \leq n \leq N,\\
      0, & \text{ in any other case.}
      \end{cases}
    \]
    Observe that $f_N \in \Hpz$, as it is finitely supported. Now, for each $N \in \N$ we have $M_p(0,f_N)=1$ and for each $n=1, 2, \dots, N+1$ we have
    \begin{align*}
  M_{p}^{p}(n,f_N)
  &= \frac{1}{\gamma(n)}\sum_{|v|=n}|f_N(v)|^{p} \\
      &=  \frac{1}{\gamma(n)}\sum_{v\in \Chi{w_{n-1}}}|f_N(v)|^{p}\\
  &=  \frac{\gamma(1,w_{n-1})}{\gamma(n)}\left|\left(
    \frac{\gamma(n)}{K(1,n-1)}\right)^{\frac{1}{p}}\right|^{p}\\
  &= 1.
\end{align*}
If $n > N+1$, then $M_p^p(n,f_N)=0$ and hence $\| f_N \|=1$.

On the other hand, as before, if $|v|=n$ but $v\neq w_n$, we have $(Bf)(v)=0$. Hence, if $0 \leq n \leq N$, we obtain
\begin{align*}
  M_p^p(n,B f_N)
  &= \frac{1}{\gamma(n)} \sum_{|v|=n} |(Bf_N)(v)|^p \\
  &= \frac{1}{\gamma(n)} |(Bf_N)(w_n)|^p \\
  &= \frac{1}{\gamma(n)} \left|\sum_{u \in \Chi{w_n}} f(u) \right|^p \\
  &= \frac{1}{\gamma(n)} \left| K(1,n) \left(\frac{\gamma(n+1)}{K(1,n)}\right)^{\frac{1}{p}} \right|^p \\
  &= \frac{(K(1,n))^{p-1}\gamma(n+1)}{\gamma(n)},
\end{align*}
while if $n > N$ we have $M_p^p(n,B f_N)=0$. This implies that
\[
  \| B f_N \|^p = \sup\left\{\frac{(K(1,n))^{p-1}\gamma(n+1)}{\gamma(n)} \, : \, 0 \leq n \leq N \right\},
\]
and, since $B$ is bounded, we get
\[
 \sup\left\{\frac{(K(1,n))^{p-1}\gamma(n+1)}{\gamma(n)} \, : \, 0 \leq n \leq N \right\} \leq \|B \|^p,
\]
for all $N \in \N$. This implies that the set
\[
  \left\{\frac{(K(1,n))^{p-1}\gamma(n+1)}{\gamma(n)} \, : \, n \in \N_0 \right\}
\]
is bounded and
\[
  \sup\left\{\frac{(K(1,n))^{p-1}\gamma(n+1)}{\gamma(n)} \, : \, n \in \N_0 \right\} \leq \| B \|^p,
\]
finishing the proof.
  \end{proof}

  Observe that since $(K(1,n))^{p-1} \geq 1$ for every $n$ and it is impossible to have  $\gamma(n+1) < \gamma(n)$ for all $n$, it follows that if $B$ is bounded, we will always have $\| B \| \geq 1$. It is also worth noticing that if $B$ is bounded and $\|B \|=1$, then
  \[
    \gamma(n+1) \leq (K(1,n))^{p-1} \gamma(n+1) \leq \gamma(n)
    \]
    for every $n \in \N_0$, which implies that $T$ is the tree where every vertex has exactly one child. Therefore, if $T$ is not this (trivial) tree, then $\| B \| > 1$. (Note that, in either case, $B$ is never an isometry.)
    
   We now show two examples, one where $B$ is bounded and one where it is not.
  
  \begin{example}
    Let $m \in \N$ and let $T$ be the tree where every vertex has $m$ children. Then $B$ is bounded on $\Hp$ and on $\Hpz$.
  \end{example}
  \begin{proof}
    Clearly, in this case, $K(1,n)=m$ and $\gamma(n)=m^n$. Hence
    \[
      \frac{(K(1,n))^{p-1}\gamma(n+1)}{\gamma(n)}=\frac{m^{p-1} m^{n+1}}{m^n}=m^p.
    \]
    Therefore $B$ is bounded and $\|B \|=m$.
    \end{proof}

    \begin{example}
      Let $T$ be the tree such that each vertex at level $n$ has exactly $n+2$ children. Then $B$ is unbounded on $\Hp$ and on $\Hpz$.
    \end{example}
    \begin{proof}
      Observe that, for each $n \in \N_0$, we have $K(1,n)=n+2$, and $\gamma(n)=(n+1)!$. Hence
          \[
      \frac{(K(1,n))^{p-1}\gamma(n+1)}{\gamma(n)}=\frac{(n+2)^{p-1} (n+2)!}{(n+1)!}=(n+2)^p,
    \]
    which is unbounded. The conclussion follows.
\end{proof}    

Later, we will need the following result to obtain the spectral radius of $B$. Observe that for $f \in \calF$ and for every $n \in \N$  we have
\[
  (B^n f)(v)=\sum_{w \in \nChi{v}{n}} f(w),
\]
for each vertex $v$.

\begin{theorem}\label{th:BmboundedHp}
  Let $T$ be a tree and let $1 \leq p < \infty$. Assume the backward shift $B$ is bounded on $\Hp$. For every $m \in \N$ we have
  \[
    \| B^m \| = \left(\sup \left\{ \frac{(K(m,n))^{p-1} \gamma(n+m)}{\gamma(n)} \, : \, n \in \N_0 \right\} \right)^{1/p}.
  \]
  \end{theorem}
  \begin{proof}
The proof follows the same steps as the proof of Theorem~\ref{th:BboundedHp} above, so we omit it.
\end{proof}

\subsection{Spectrum of $B$}\label{subsec_specB} We will now study the spectrum of the backward shift $B$. Although we restrict ourselves to a certain type of tree, as will be seen below, we start with an observation about a particular type of eigenfunctions in $\calF$ for $B$ in arbitrary trees.

\begin{lemma}\label{le:eigenfunctions}
Let $T$ be a tree and let $\lambda \in \C$. If we define $f\in\calF$ as
\[
    f(v):=\begin{cases}
    \frac{\lambda^{|v|}}{\gamma(1,\parent{v}) \gamma(1,\nparent{v}{2}) \dots \gamma(1,\nparent{v}{|v|-1}) \gamma(1,\root)},& \text{ if } v\neq\root,\\
    1, & \text{ if } v = \root,
    \end{cases}
  \]
  then $B f = \lambda f$.
\end{lemma}
\begin{proof}
  If $v$ is a vertex then
\begin{eqnarray*}
    Bf(v)
    &=&\sum_{w\in \Chi{v}}f(w)\\
    &=&\sum_{w\in \Chi{v}}\frac{\lambda^{|w|}}{\gamma(1,\parent{w})\gamma(1,\nparent{w}{2})\dots\gamma(1,\nparent{w}{|w|-1})\gamma(1,o)}\\
    &=&\sum_{w\in \Chi{v}}\frac{\lambda^{|v|+1}}{\gamma(1,v)\gamma(1,\parent{v})\dots\gamma(1,\nparent{v}{|v|-1})\gamma(1,o)}\\ &=&\gamma(1,v)\frac{\lambda^{|v|+1}}{\gamma(1,v)\gamma(1,\parent{v})\dots\gamma(1, \nparent{v}{|v|-1})\gamma(1,o)}\\  &=&\lambda\frac{\lambda^{|v|}}{\gamma(1,\parent{v})\dots\gamma(1,\nparent{v}{|v|-1})\gamma(1,o)}\\
    &=&\lambda f(v).
\end{eqnarray*}
This proves the desired equality.
\end{proof}

Based on the above result, it can be easily checked that $\conj{\D}$ is always contained in the set of eigenvalues of $B$ as an operator on $\Hp$. Indeed, just observe that if $|\lambda|\leq 1$, the function $f$ in Lemma~\ref{le:eigenfunctions} is bounded. Also, observe that since $|f(v)|\leq |\lambda|^{|v|}$ and hence $M_p(n,f) \leq |\lambda|^n$, in the case  $|\lambda|<1$ we have that $M_p(n,f)  \to 0$ as $n \to \infty$. Hence $\D$ is contained in the set of eigenvalues of $B$ as an operator on $\Hpz$. This is summarized in the following proposition.

\begin{proposition}\label{prop:disk_spec}
Let $T$ be a tree. Then as an operator on $\Hp$, we have $\conj{\D} \subseteq \sigma_p(B)$ and as an operator on $\Hpz$ we have $\D \subseteq \sigma_p(B)$.
\end{proposition}

In case the tree has the property that all vertices in the same level have the same number of children, we can give better lower and upper estimates for the spectrum of $B$ on $\Hpz$. We restrict ourselves to this type of trees for the remainder of this study of the spectrum.

\begin{theorem}\label{th:estimate_spectrum}
  Let $T$ be a tree and assume there exists a bounded sequence of positive integers $(s_n)$ such that $\gamma(1,v)=s_{|v|+1}$ for all vertices $v$. Then $B$ is a bounded operator and $\| B \|= \sup_{n\in \N} s_n$. Furthermore,
  \[
    \left\{ \lambda \in \C \, : \, \left(\frac{\lambda^{n}}{s_1 s_2 \dots s_n} \right) \text{ is bounded } \right\} \subseteq \sigma_p(B),
  \]
  where $\sigma_p(B)$ denotes the set of eigenvalues of $B$ as an operator on $\Hp$. Also,  the spectral radius of $B$ as an operator on $\Hp$ equals
  \[
    r(B)=\lim_{m\to \infty} \left( \sup_{n \in \N_0} s_{n+1} s_{n+2} \cdots s_{n+m} \right)^{\frac{1}{m}}.
    \]
  \end{theorem}
  \begin{proof}
    The hypothesis implies that for each $n \in \N_0$ we have $K(1,n)=s_{n+1}$ and that $\gamma(n)=s_1 s_2 \dots s_n$. Hence
    \[
      \frac{(K(1,n))^{p-1} \gamma(n+1)}{\gamma(n)} = \frac{(s_{n+1})^{p-1} s_1 s_2 \dots s_{n+1}}{s_1 s_2 \dots s_n} = s_{n+1}^p,
    \]
    which is bounded by hypothesis. Hence $B$ is a bounded operator, with norm $\|B \|=\sup_{n\in \N} s_n$.

    Now, in this case, the functions on Lemma~\ref{le:eigenfunctions} can be written as
    \[
      f(v)=\begin{cases}
        \frac{\lambda^{|v|}}{s_{|\!v\!|}  \, \cdots \, s_2 \,  s_1},& \text{ if } v\neq \root,\\
        1,& \text{ if } v = \root,
        \end{cases}
      \]
      for each vertex $v$. Let us verify that $f \in \Hp$: indeed, for each $n \in \N_0$ we have
\begin{align*}
        M_{p}^{p}(n,f)
  &=\frac{1}{\gamma(n)}\sum_{|v|=n}|f(v)|^{p} \\
  &=\frac{1}{\gamma(n)}\sum_{|v|=n}\left| \frac{\lambda^{n}}{s_{1}s_{2} \cdot\cdot\cdot s_{n}}  \right|^{p} \\
  &=\left(\frac{|\lambda|^{n}}{s_{1}s_{2} \cdot\cdot\cdot s_{n}}\right)^{p}.
\end{align*}
The hypothesis implies that $M_p(n,f)$ is bounded and therefore $f \in \Hp$. This proves the first part of the theorem.

For the second part, observe that $K(m,n)=s_{n+1} s_{n+2} \dots s_{n+m}$. Hence Theorem~\ref{th:BmboundedHp} implies that
\begin{align*}
  \| B^m \|^p
  &= \sup_{n\in \N_0} \frac{(K(m,n))^{p-1}\gamma(n+m)}{\gamma(n)}  \\
  &= \sup_{n\in \N_0} \frac{(s_{n+1} s_{n+2} \dots s_{n+m})^{p-1} s_1 s_2 \dots s_{n+m}}{s_1 s_2 \dots s_n} \\
  &= \sup_{n\in \N_0} (s_{n+1} s_{n+2} \dots s_{n+m})^p
\end{align*}
and therefore
\[
  \| B^m \| = \sup_{n\in \N_0} s_{n+1} s_{n+2} \dots s_{n+m}.
\]
But the spectral radius theorem guarantees that
\[
  \lim_{m\to \infty }\| B^m \|^{\frac{1}{m}} = \left(\sup_{n\in \N_0} s_{n+1} s_{n+2} \dots s_{n+m}\right)^{\frac{1}{m}}
\]
exists and equals $r(B)$, finishing the theorem.
\end{proof}

Observe that the proof above also shows that
  \[
    \left\{ \lambda \in \C \, : \, \lim_{n\to \infty } \frac{\lambda^{n}}{s_1 s_2 \dots s_n} =0  \right\} \subseteq \sigma_p(B),
  \]
  where now $\sigma_p(B)$ denotes the set of eigenvalues of $B$ as an operator on $\Hpz$. The spectral radius of $B$ as an operator on $\Hpz$ clearly satisfies the same formula as above.

The following corollary improves, when the hypothesis of the previous theorem hold, the results of Proposition~\ref{prop:disk_spec}.
  
  \begin{corollary}\label{cor:disk_spec_root}
    Let $T$ be a tree and assume there exists a bounded sequence of positive integers $(s_n)$ such that $\gamma(1,v)=s_{|v|+1}$ for all vertices $v$. If $\ds t:=\liminf_{n\to \infty} (s_1 s_2 \dots s_n)^{\frac{1}{n}}$, then $t \D \subseteq \sigma_p(B)$, both as an operator on $\Hp$ and on $\Hpz$.
  \end{corollary}
  \begin{proof}
    By the root test, if $\lambda \in t \D$, the sequence
    \[
      \left(\frac{\lambda^{n}}{s_1 s_2 \dots s_n} \right)
    \]
    tends to zero since
    \[
      \limsup_{n\to \infty} \left(\frac{|\lambda|^{n}}{s_1 s_2 \dots s_n}\right)^{\frac{1}{n}}
        =\limsup_{n\to \infty} \frac{|\lambda|}{(s_1 s_2 \dots s_n)^{\frac{1}{n}}}= \frac{|\lambda|}{\ds\liminf_{n \to \infty } (s_1 s_2 \dots s_n)^{\frac{1}{n}} }=\frac{|\lambda|}{t} <1.
    \]
Theorem~\ref{th:estimate_spectrum} then proves the corollary.
    \end{proof}

    Since $\liminf_{n\to \infty} (s_1 s_2 \dots s_n)^{\frac{1}{n}} \leq \sup_{n \in \N} s_n$, the above result allows us to obtain, in case we have an equality, a complete description of the spectrum of $B$.
    
  \begin{corollary}\label{cor:inf_sup}
    Let $T$ be a tree and assume there exists a bounded sequence of positive integers $(s_n)$ such that $\gamma(1,v)=s_{|v|+1}$ for all vertices $v$. If $\ds t:=\liminf_{n\to \infty} (s_1 s_2 \dots s_n)^{\frac{1}{n}}=\sup_{n \in \N} s_n$, then $\sigma(B)=\conj{t \D}$, both as an operator on $\Hp$ and on $\Hpz$.
  \end{corollary}

We will now explore some other cases where we can completely describe the spectrum of $B$, both as an operator on $\Hp$ and as an operator on $\Hpz$.

\begin{example}
  Let $T$ be a tree where every vertex has $q$ children. Then $\sigma(B)=\conj{q \D}$, both as an operator on $\Hp$ and on $\Hpz$.
\end{example}
\begin{proof}
  This follows easily in the case where we consider $B$ as an operator on $\Hp$ by applying the above theorem and noticing that $s_n=q$ for all $n \in \N$.

  In this case,
  \[
  \liminf_{n\to \infty} (s_1 s_2 \dots s_n)^{\frac{1}{n}}= q = \sup_{n \in \N} s_n,
\]
and hence by Corollary~\ref{cor:inf_sup}, we have $\sigma(B)=\conj{q \D}$.
\end{proof}

\begin{example}
  Let $T$ be a tree and let $q_1, q_2, \dots, q_m$ be positive integers such that, for all vertices $v$, we have $\gamma(1,v)=q_{r+1}$, where $r$ is an integer such that $|v|=m p + r$, with $0\leq r < m$ for some integer $p$. Then $\sigma(B)$ is the closed disk centered at the origin with radius $(q_1 \cdot q_2 \cdots q_m)^{\frac{1}{m}}$.
  \end{example}
  \begin{proof}
    Let $A:=(q_1 \cdot q_2 \cdots q_m)^{\frac{1}{m}}$ and make $s_n=\gamma(1,v)$, where $v$ is any vertex with $|v|=n-1$. If $p$ and $n$ are integers such that $n=m p +r$ with $0 \leq r < m$, then
    \[
      (s_1 s_2 \dots s_n)^{\frac{1}{n}} =  \left( (q_1 q_2 \cdots q_m)^{p} (q_1 q_2 \cdots q_{r}) \right)^{\frac{1}{n}}= \left( A^{m p} (q_1 q_2 \cdots q_{r}) \right)^{\frac{1}{n}}= A \left( A^{-\frac{r}{n}}  (q_1 q_2 \cdots q_{r})^{\frac{1}{n}}\right).
    \]
    But observe that
    \[
      1 \leq (q_1 q_2 \cdots q_{r})^{\frac{1}{n}} \leq (\max\{q_1, q_2, \dots q_m\})^{\frac{m}{n}},
    \]
    and hence $\ds \lim_{n \to \infty} (q_1 q_2 \cdots q_{r})^{\frac{1}{n}} =1$. Also observe that
    \[
      A ^{-\frac{m}{n} }\leq A^{-\frac{r}{n}} \leq 1,
    \]
    since $A\geq 1$ and $ r < m$, and hence  $\ds \lim_{n \to \infty} A^{-\frac{r}{n}}=1$. Therefore,
    \[
      \liminf_{n\to \infty} (s_1 s_2 \dots s_n)^{\frac{1}{n}} = A,
    \]
    which shows that $A \D \subseteq \sigma_p(B)$.

    Observe also that, since we know
    \[
      r(B)=\lim_{j\to \infty} \left(\sup_{n\in \N_0} s_{n+1} s_{n+2} \dots s_{n+j}\right)^{\frac{1}{j}}
    \]
    exists, we may consider the limit of a subsequence. So make $j=km$ and observe that
    \[
      \left(\sup_{n\in \N_0} s_{n+1} s_{n+2} \dots s_{n+km}\right)^{\frac{1}{km}} = \left(\sup_{n\in \N} \  (q_1 q_2 q_3 \dots q_m)^k \right)^{\frac{1}{km}} = \left(q_1 q_2 q_3 \dots q_m \right)^{\frac{1}{m}} = A.
    \]
    Hence, $r(B)=A$. This proves the result.
    \end{proof}

    In the example above, the spectral radius is the geometric mean of the (finite) sequence of degrees. In the example below we see that, in general, our results are not sharp enough to find the spectrum of $B$.
    
    \begin{example}
      Let $T$ be the tree such that $\gamma(1,v)=s_{|v|+1}$, where the sequence $(s_n)$ is given by the list $(2, 3, 2, 2, 3, 3, 2, 2, 2, 3, 3, 3, \dots)$; that is
      \[
        s_n:=\begin{cases}
          3, & \text{ if } k^2 +1 \leq n \leq  k^2+k \text{ for some } k \in \N, \\
          2, & \text{ in any other case.}
        \end{cases}
      \]
      Then $r(B)=\lim_{m\to \infty} \left(  \sup_{n \in \N_0} s_{n+1} s_{n+2} \cdots s_{n+m} \right)^{1/m}= 3$ and $\liminf_{n\to \infty} (s_1 s_2 \cdots s_n)^{\frac{1}{n}} \leq \sqrt{6}$.
\end{example}
\begin{proof}
  Observe that, for each fixed $m \in \N$, we have
  \[
    s_{m^2+1} s_{m^2+2} \cdots s_{m^2+m} = 3^m,
  \]
  and hence
  \[
    \sup_{n \in \N_0} s_{n+1} s_{n+2} \cdots s_{n+m}=3^m.
  \]
  This in turn implies that $r(B)=\lim_{m\to \infty} \left(  \sup_{n \in \N_0} s_{n+1} s_{n+2} \cdots s_{n+m} \right)^{1/m} = 3$.

  On the other hand, if we consider the sequence given by $n_k=k(k+1)$, then
  \[
    s_1 s_2 \cdots s_{n_k} = 2^{\frac{k(k+1)}{2}} 3^{\frac{k(k+1)}{2}}
  \]
  and hence
  \[
    \lim_{k\to \infty} \left( s_1 s_2 \cdots s_{n_k} \right)^{\frac{1}{n_k}} = 2^{\frac{1}{2}} 3^{\frac{1}{2}}.
  \]
  But this implies that
  \[
    \liminf_{n\to \infty} (s_1 s_2 \cdots s_n)^{\frac{1}{n}} \leq \sqrt{6}. \qedhere
  \]
\end{proof}

Is it possible to obtain a complete characterization of the spectrum of $B$ in terms of the sequence $(s_n)$? We leave this question open for future research.


\subsection{Hypercyclicity of $B$}\label{subsec_hyperB}
We now study the hyperciclicity of $B$ on $\Hpz$. First we make an easy observation.
    \begin{proposition}
      Let $T$ be a tree. If $T$ has a leaf, then $B$ is not hypercyclic.
    \end{proposition}
    \begin{proof}
      Let $v$ be a leaf. Observe that for any $f \in \Hpz$, we have
      \[
        (B^n f)(v)=\sum_{w \in \nChi{v}{n}} f(w) = 0,
      \]
      since $\nChi{v}{n}=\varnothing$ for all $n \in \N$. Set $s=|v|$. If $B$ were hypercyclic, with hypercyclic vector $f \in \Hpz$, there would exist $N \in \N$ with
      \[
        \| B^N f - \chi_{\{v\}} \| < \frac{1}{2 (\gamma(s))^{1/p}}.
      \]
      But since
      \[
        \frac{1}{(\gamma(s))^{1/p}} | (B^N f)(v) - \chi_{\{v\}}(v)| \leq \left( \frac{1}{\gamma(s)} \sum_{ |w|=s}| (B^N f)(w) - \chi_{\{v\}}(w)|^p \right)^{1/p}
       \leq  \| B^N f - \chi_{\{v\}} \|,
     \]
    we obtain
     \[
 \frac{1}{(\gamma(s))^{1/p}} | 0 -1 | =   \frac{1}{(\gamma(s))^{1/p}} | (B^N f)(v) - \chi_{\{v\}}(v)| < \frac{1}{2 (\gamma(s))^{1/p}},
     \]
which is a contradcition. Hence $B$ cannot be hypercyclic.       
      \end{proof}

      As it turns out, hyperciclicity of $B$, in the leafless case, is equivalent to the number of vertices in each level going to infinity. Compare to what happens in other spaces: for example, in the Lipschitz space of the tree (see \cite{MaRi}) $B$ cannot be hypercyclic if there are free ends; here they are allowed.
      
    \begin{theorem}
      Let $T$ be a leafless tree. Then $B$ is hypercyclic on $\Hpz$ if and only if $\ds \lim_{n\to \infty}\gamma(n)=\infty$.
    \end{theorem}
      \begin{proof}
        Assume $B$ is hypercyclic. Let $\chi_{\root}$ be the characteristic function of the root $\root$. Choose a positive number $M \in \R$ and let $\epsilon:=\frac{1}{1+M}$. 

        Let $f$ be a hypercyclic vector for $B$ with $\| f \| < \epsilon$. By hypercyclicity of $B$, there exists $N \in \N$ such that
        \[
          \| B^N f - \chi_{\root} \| < \epsilon.
        \]
        We then have
        \[
          | (B^N f)(\root) - 1|^p
          \leq \sup_{n \in \N_0} \frac{1}{\gamma(n)} \sum_{|v|=n} |(B^N f)(v) - \chi_{\root}(v) |^p 
          = \| B^N f - \chi_{\root} \|^p 
          < \epsilon^p.
        \]
This implies that
          \[
            \left| \left(\sum_{|v|=N} f(v) \right) - 1\right| < \epsilon,
          \]
          and hence
        \[
          1 - \epsilon <  \left| \sum_{|v|=N} f(v) \right| \leq \sum_{|v|=N} |f(v)|.
        \]
        Dividing by $\gamma(N)$ and using Jensen's inequality we obtain
        \[
          \frac{1}{\gamma(N)} (1-\epsilon) <  \frac{1}{\gamma(N)} \sum_{|v|=N} |f(v)| \leq  \left( \frac{1}{\gamma(N)} \sum_{|v|=N} |f(v)|^p \right)^{1/p} \leq \| f \|.
        \]
        Hence
        \[
          \frac{1}{\gamma(N)} (1-\epsilon) < \| f \| < \epsilon.
        \]
        Since $\epsilon=\frac{1}{1+M}$ this implies that
        \begin{equation*}\label{eq:limit}
          M < \gamma(N).
        \end{equation*}
        Since $T$ is leafless, it follows that $\gamma(\cdot)$ is nondecreasing. Hence, since we have shown that for each positive $M \in \R$ there exists $N \in \N$, with $M < \gamma(N)$, this implies that $\ds \lim_{n\to \infty} \gamma(n)=\infty$.

        Now, assume that $\ds \lim_{n\to \infty} \gamma(n)=\infty$. We will prove the hypercyclicity of $B$ by using the Kitai-Gethner-Shapiro criterion for hyperciclicity (see, for example, \cite[pp.74--75]{GrPe}).

        First, let $X=\{ f \in \Hpz \, : \, f \text{ has finite support} \}$. It is straightforward to check that $X$ is dense in $\Hpz$.
        \begin{itemize}
        \item Let $f \in X$ and let $N \in \N$ be such that $f(w)=0$ for all $w$ with $|w|>N$. Then, if $v$ is a vertex and $n > N$ we have
          \[
            (B^n f)(v)=\sum_{w \in \nChi{v}{n}} f(w) =0;
          \]
          i.e., $B^n f = 0$ for $n \geq N$. Hence $B^n f \to 0$ as $n \to \infty$, as required.
        \item For each $g \in X$, and every $n \in \N$, we define a function $T_n g \in \calF$ as 
          \[
            (T_n g)(v)=\begin{cases}
              \frac{g(\nparent{v}{n})}{\gamma(\nparent{v}{n},n)}, & \text{ if } v \text{ has an $n$-parent}, \\
              0, & \text{ if } v \text{ does not have an $n$-parent.}
            \end{cases}
          \]
          Since $g$ is of finite support, then $T_n g$ is also of finite support and hence is in $\Hpz$. Let $m \in \N_0$. If $m < n$, then clearly $M_p(m,T_ng)=0$. If $m \geq n$, then
            \begin{align*}
              M_p(m, T_n g)
              &=\frac{1}{\gamma(m)} \sum_{|v|=m} |(T_ng)(v)|^p \\
              &=\frac{1}{\gamma(m)} \sum_{|v|=m} \left|\frac{g(\nparent{v}{n})}{\gamma(\nparent{v}{n},n)}\right|^p \\
              &=\frac{1}{\gamma(m)} \sum_{|v|=m-n} \left|\frac{g(v)}{\gamma(v,n)} \right|^p \gamma(v,n) \\
              &=\frac{1}{\gamma(m)} \sum_{|v|=m-n} \frac{1}{(\gamma(v,n))^{p-1}} |g(v)|^p.
            \end{align*}
            Since the tree is leafless, $\gamma(v,n) \geq 1$ and hence
            \begin{align*}
              M_p^p(m,T_n g)
              &\leq \frac{1}{\gamma(m)} \sum_{|v|=m-n} |g(v)|^p \\
              &= \frac{\gamma(m-n)}{\gamma(m)} M_p^p(m-n,g).
            \end{align*}

The inequality above implies that
            \[
             \| T_n g \|^p =\sup_{m \geq n} M_p^p(m, T_n g) \leq \sup_{m \geq n} \frac{\gamma(m-n)}{\gamma(m)} M_p^p(m-n,g).
            \]

            Now, let $N$ be such that $g(v)=0$ if $|v|> N$. We then have
\begin{align*}
  \| T_n g \|^p &\leq\sup_{0 \leq m-n \leq N} \frac{\gamma(m-n)}{\gamma(m)} M_p^p(m-n,g)\\
              &=\sup_{0 \leq s \leq N} \frac{\gamma(s)}{\gamma(s+n)} M_p^p(s,g)\\
  & \leq \max_{0 \leq s \leq N} \frac{\gamma(s)}{\gamma(s+n)} \| g \|^p.
            \end{align*}
            Since $\gamma(s+n) \to \infty$  as $n \to \infty$ for each $s$ with $0 \leq s \leq N$, we have $ \| T_n g \| \to 0$ as $n \to \infty$.

          \item Lastly, observe that
            \begin{align*}
              (B^n T_n g) (v)
              &= \sum_{w \in \nChi{v}{n}} (T_n g)(w) \\
              &= \sum_{w \in \nChi{v}{n}}  \frac{g(\nparent{w}{n})}{\gamma(\nparent{w}{n},n)} \\
              &=\sum_{w \in \nChi{v}{n}} \frac{g(v)}{\gamma(v,n)}\\
              &=g(v).
            \end{align*}
            Hence, $B^n T_n g= g$ and the last condition of the hypercyclicity criterion is satisfied.            
          \end{itemize}
        This concludes the proof of the hypercyclicity of $B$.
      \end{proof}


\bibliographystyle{amsplain}

\bibliography{trees_hardy_shift}{}

    \end{document}